\documentclass[12pt,leqno]{amsart}

\usepackage{enumerate}

\usepackage{eucal}
\usepackage{euler}
\usepackage{times}

\swapnumbers

\theoremstyle{plain}
    \newtheorem{theorem}{Theorem}
    \newtheorem{lemma}[theorem]{Lemma}

 \newtheorem{definition-theorem}[theorem]{Definition/Theorem}

\theoremstyle{definition}

\theoremstyle{remark}

\newcounter{ictr}
\newenvironment{ilist}{\begin{list}
                         {\textup{(\arabic{ictr})}}
                         {\usecounter{ictr}
                          \setlength{\leftmargin}{0.6truein}
                          \setlength{\itemsep}{0.00truein}
                          \setlength{\labelwidth}{0.3truein}}}
                      {\end{list}}

\newcounter{kctr}
\newenvironment{klist}{\begin{list}
                         {\textup{(\alph{kctr})}}
                         {\usecounter{kctr}
                          \setlength{\leftmargin}{0.6truein}
                          \setlength{\itemsep}{0.00truein}
                          \setlength{\labelwidth}{0.3truein}}}
                      {\end{list}}

\newcommand{\Hawaii}{Hawai\kern.05em`\kern.05em\relax i}
\newcommand{\Manoa}{M\=anoa}

\newcommand{\R}{\mathbb{R}}

\newcommand{\e}{\mathcal{E}}

\newcommand{\eps}{\varepsilon}
    
\DeclareMathOperator{\supp}{supp}

\begin{document}   

\title[Uniform embeddings into reflexive Banach space]{
Uniform embeddings of bounded geometry spaces into reflexive Banach space
}


\author{Nathanial Brown}
\address{Department of Mathematics, Pennsylvania State University,
  University Park, PA 16802}
\email{nbrown@math.psu.edu}
\author{Erik Guentner}
\address{University of \Hawaii, \Manoa, Department of Mathematics, 2565
  McCarthy Mall, Honolulu, HI 96822-2273}
\email{erik@math.hawaii.edu}

\thanks{The authors were partially supported by  grants from the
U.S. National Science Foundation. }

\begin{abstract}
  We show that every metric space with bounded geometry uniformly embeds
  into a direct sum of $l^p ({\mathbb N})$ spaces ($p$'s going off to
  infinity).  In particular, every sequence of expanding graphs
  uniformly embeds into such a reflexive Banach space even though no
  such sequence uniformly embeds into a fixed $l^p ({\mathbb N})$ space.
  In the case of discrete groups we prove the analogue of
  a-$T$-menability -- the existence of a metrically proper affine
  isometric action on a direct sum of $l^p ({\mathbb N})$ spaces.
\end{abstract}

\maketitle

\section*{Introduction}
\label{sec:introduction}

Gromov introduced the notion of uniform embeddability and later
suggested that uniform embeddability of a discrete group into a Hilbert
space, or even into a uniformly convex Banach space, might be relevant
for applications to the Novikov conjecture \cite{gromov93,NCITR}.  Yu
subsequently proved the Coarse Baum-Connes and Novikov conjectures for
groups which are uniformly embeddable into a Hilbert space
\cite{yu00,skandalis-tu-yu02}.  This lead to the verification of the
Novikov conjecture for large classes of discrete groups, including, for
example, linear groups \cite{guentner-higson-weinberger-unpub} and, more
generally, exact groups \cite{ozawa00,guentner-kaminker02}.

The question of whether every discrete group might be uniformly
embeddable in Hilbert space remained open for some time.  Gromov ended
the speculation by announcing the construction of a discrete group which
cannot be uniformly embedded into Hilbert space
\cite{gromov-random,gromov00}.  Nevertheless, the Novikov conjecture
remains valid for the groups Gromov constructs.

Two natural questions arise: (1) Can one prove the Novikov
conjecture for discrete groups which uniformly embed into other nice
Banach spaces?  (2) Does every discrete group embed into a nice Banach
space?  How `nice' the range Banach space needs to be (for question (1))
is not yet clear, although quite recently Kasparov and Yu have achieved
exciting results in this direction.  In this note we show that
every metric space with bounded geometry can be uniformly embedded into a
fairly nice Banach space.  More precisely, our main result is the
following:

\begin{theorem}
\label{thm:bddgeomcase} Let $X$ be a metric space with bounded
geometry.  There exists a sequence of positive real numbers
$\{p_n\}$ and a uniform embedding of $X$ into the $l^2$-direct sum
$\oplus l^{p_n}({\mathbb N})$.
\end{theorem}

The Banach space in the statement is reflexive.  Thus, every bounded
geometry metric space can be uniformly embedded into a reflexive Banach
space.  In particular, every sequence of expanding graphs can be
uniformly embedded into a reflexive Banach space, answering a question
pondered by several experts.  As observed by Gromov, a sequence of
expanding graphs cannot be uniformly embedded into Hilbert space
\cite{gromov00} (see \cite{guentner-higson-cime} for a proof); a less
well-known observation by Roe is that such a sequence cannot be
uniformly embedded into $l^p({\mathbb N})$ for any finite $p$ (cf.\ 
\cite{matousek97}).  Hence the result above seems to give a uniform
embedding into the simplest kind of Banach space that is capable of
containing all bounded geometry metric spaces.

In the case of discrete groups we prove the following stronger result.

\begin{theorem}
\label{thm:groupcase} 
Let $\Gamma$ be a countable discrete group.  There exists a sequence of
positive real numbers $\{p_n\}$ and a metrically proper affine isometric
action of $\Gamma$ on the $l^2$-direct sum 
$\oplus l^{p_n}({\mathbb N})$.  
In particular, there exists a uniform embedding 
$\Phi : \Gamma \to \oplus l^{p_n}({\mathbb N})$ such that 
$\| \Phi(s) - \Phi(t) \| = \| \Phi(t^{-1}s) \|$, for all $s$,
$t\in\Gamma$. 
\end{theorem}

This theorem is a version of a-$T$-menability involving a reflexive Banach
space in place of a Hilbert space.  The remarks above apply equally to
the groups constructed by Gromov: they cannot be uniformly embedded into,
and in particular do not admit metrically proper affine isometric actions
on, $l^p ({\mathbb N})$ for any finite $p$.  Thus, again, the result is in
some sense optimal.

The proofs of the two results above have certain similarities though the
actual constructions involved are slightly different.  In the next
section we collect several facts needed for both proofs and set our
notation.  In the remaining two sections we construct the required
embedding and affine action, proving the theorems.

\section{Notation, Definitions and Basic Facts}
\label{sec:notat-defin-basic}

For a countable set $S$ let $l^p (S)$ denote the usual sequence space
and $\| \cdot \|_p$ the corresponding norm (we allow $1 \leq p \leq
\infty$).  For a sequence ${\mathcal E}_n$ of Banach spaces let 
$\oplus {\mathcal E}_n$ denote the $l^2$-direct sum defined by 
$\oplus {\mathcal E}_n = \{ (z_n): \sum \|z_n \|^2 < \infty \}$, with
the obvious norm.   Since the dual space of an $l^2$-direct sum is the
$l^2$-direct sum of the dual spaces we see that if every 
${\mathcal E}_n$ is reflexive then so is $\oplus {\mathcal E}_n$.  In
particular, for every sequence of real numbers $\{p_n\}$, 
$\oplus l^{p_n}({\mathbb N})$ is a reflexive
Banach space.

A metric space $X$ is \emph{locally finite} if every ball is finite; it
has \emph{bounded geometry} if for all $R>0$ there exists $C>0$ such
that every $R$-ball in $X$ contains at most $C$ points.  With these
conventions a bounded geometry metric space is necessarily discrete and
countable.

A function $\Phi:X\to \e$ from a metric space $X$ to a Banach space $\e$
is a \emph{uniform embedding} if there exist non-decreasing proper
functions $\rho_\pm:[0,\infty) \to \R$ such that
\begin{equation*}
  \rho_-(d(x,y)) \leq \| \Phi(x)-\Phi(y) \| \leq \rho_+(d(x,y)), \quad
        \text{for all $x$, $y\in X$},
\end{equation*}
where $\|\cdot\|$ denotes the norm in $\e$.


We require the following elementary fact.

\begin{lemma}
\label{lem:conv} Let $S$ be a set, and let $\alpha\geq 0$ 
and $\beta\geq 1$.  Then $\|f\|_p\to \|f\|_\infty$ as $p\to\infty$
uniformly on sets of the form
\begin{equation*}
  \mathcal S_{\alpha,\beta} =
  \{\, f\in l^\infty(S) \,\colon\, \text{$\|f\|_\infty\leq \alpha$
              and $\#\supp(f) \leq \beta$} \,\}, 
\end{equation*} where $\#\supp(f)$ denotes the cardinality of the
support of $f$.
\end{lemma}

\begin{proof} The following inequalities imply the lemma:
\begin{align*}
   \| f \|_{\infty} &\leq \| f\|_p  = \left\{ \sum_{s \in supp(f)} 
                                        |f(s)|^p \right\}^{1/p} \\ 
            &\leq \| f \|_{\infty} \left\{ \sum_{s 
                        \in supp(f)} 1 \right\}^{1/p}  
                     \leq \| f \|_{\infty} \beta^{1/p}. \quad\qed
\end{align*}
\renewcommand{\qed}{}\end{proof}

\section{Bounded Geometry Case}
\label{sec:bound-geom-case}

In this section we prove Theorem~\ref{thm:bddgeomcase}.  The main
observation is that every locally finite metric space enjoys a weak
analogue of Yu's property A \cite{yu00}.  More precisely we have the
following lemma.

\begin{lemma}
\label{lem:bddgeomlem}
Let $X$ be a locally finite metric space. Then there exists a
sequence of mappings $\phi^{n} : X \to l^{\infty}(X)$, $x \mapsto
\phi^{n}_x$ such that 
\begin{ilist}
  \item $\| \phi^{n}_x\|_\infty = 1$, for all $x \in X$ and 
        $n \in {\mathbb N}$;
  \item $\supp (\phi^{n}_x) \subset B_n(x)$, where $B_n(x)$ is the 
        ball of radius $n$ about $x$;
  \item $\| \phi^{n}_x - \phi^{n}_y \|_\infty\leq d(x,y)/n$, for all
            $x$, $y\in X$. 
\end{ilist}
\end{lemma}

\begin{proof} Define functions $\phi^n : X \to l^{\infty}(X)$ as follows: 
\begin{equation*}
  \phi^n_x(y) = \begin{cases} 1-\frac{d(x,y)}{n}, & d(x,y)\leq n \\
                                0, & d(x,y)\geq n. \end{cases} 
\end{equation*}
Clearly conditions (1) and (2) are satisfied.  Condition (3) 
is easily check by considering cases.

\noindent\underline{{\it Case $z\notin B_n(x)\cup B_n(y)$\/}\,:\,} 
$\phi^n_x(z)=0$ and $\phi^n_y(z)=0$.

\noindent\underline{{\it Case $z\in B_n(x)\cap B_n(y)$\/}\,:\,} 
By the triangle inequality
\begin{equation*}
  |\phi^n_x(z) - \phi^n_y(z)|= |d(x,z)-d(y,z)|/n \leq d(x,y)/n.
\end{equation*}

\noindent\underline{{\it Case $z\in B_n(x) - B_n(y)$\/}\,:\,}
Again by the triangle inequality
\begin{equation*}
    |\phi^n_x(z) - \phi^n_y(z)| = (n-d(x,z))/n \leq 
 \left(d(y,z)-d(x,z)\right)/n \leq d(x,y)/n.
\end{equation*}
The remaining case is similar.
\end{proof}

The conditions (1)--(3) in the lemma are indeed a weak form of
Property~A.  For a bounded geometry space $X$ Tu \cite{tu01} showed that
Property A is equivalent to the existence of a sequence of functions
$\psi^n:X\to l^2(X)$ satisfying the following conditions:
\begin{klist}
  \item $\|\psi^n_x\|_2 = 1$, $\forall$ $x\in X$ and 
                 $\forall$ $n\in {\mathbb N}$;
  \item $\forall$ $n\in {\mathbb N}$ $\exists$ $R>0$ such that
          $\supp(\psi^n_x)\subset B_R(x)$;
  \item $\forall$ $C>0$ $\|\psi_x^n-\psi_y^n\|_2 \to 0$ uniformly for
           $x$, $y\in X$ satisfying $d(x,y)\leq C$.
\end{klist}
We note that conditions (2) and (3) in the lemma could be replaced by
direct analogues of conditions (b) and (c) (for example, replace the
precise estimate in (3) by a statement involving only uniform
convergence as in (c)) at the cost of complicating some of the arguments
below.

\begin{proof}[Proof of Theorem~\ref{thm:bddgeomcase}.]
Since $X$ is assumed to be of bounded geometry the cardinality of the
supports of the functions $\phi^{n}_x$ (defined as in the previous
lemma) are uniformly bounded above for each $n$ ({\em independent of
  x}).  Hence for each integer $n$ we can, by Lemma~\ref{lem:conv},
find a real number $p_n$ such that
\begin{equation*}
  \| \phi^n_x - \phi^n_y \|_{p_n} \leq 
        \| \phi^n_x - \phi^n_y \|_{\infty} + 1/n,
\end{equation*}
for all $x,y \in X$.  Combined with condition (3) of the lemma this
implies, in particular, that 
\begin{equation*}
  \| \phi^n_x - \phi^n_y \|_{p_n} \leq \frac{d(x,y) + 1}{n},
\end{equation*}
for all $n\in {\mathbb N}$ and  $x,y \in X$.

Now fix a base point $x_0 \in X$ and define, for all $x\in X$,
\begin{equation*}
  \Phi(x) = \oplus \phi^n_x - \phi^n_{x_0}.
\end{equation*}
A simple calculation reveals that
\begin{equation*}
   \| \Phi(x) - \Phi(y) \|^2 = 
       \sum_{n \in {\mathbb N}} \|\phi^{n}_x - \phi^{n}_y \|^2_{p_{n}} 
         \leq (d(x,y) + 1)^2 \sum_{n \in {\mathbb N}} \frac{1}{n^2}.
\end{equation*}
Setting $y=x_0$ and observing that $\Phi(x_0)=0$ we conclude that indeed
$\Phi(x)\in \oplus l^{p_n}(X)$ for all $x\in X$.  It remains to see that
$\Phi:X\to \oplus l^{p_n}(X)$ is a uniform embedding.  Consulting the
previous inequality once again we see that we may define $\rho_{+} (r) =
(\sum_{n \in {\mathbb N}} \frac{1}{n^2})^{1/2} (r + 1)$.



We define a lower distortion function by 
\begin{equation*}
  \rho_{-} (r) = \inf_{d(x,y)\geq r} \|\Phi(x)-\Phi(y)\|,
\end{equation*}
and are left to check that $\rho_{-} (r) \to \infty$ as $r \to \infty$.
In other words, for every integer $R > 0$ we must find an $r_0$ such
that $\| \Phi(x)-\Phi(y)\| > R$ for all $x,y \in X$ satisfying $d(x,y) >
r_0$.  This is easy to do: since each $\phi^{n}_x$ is supported in the
ball of radius $n$ about $x$ if $d(x,y) > 2R^2$ it follows that
\begin{equation*}
   \| \Phi(x)-\Phi(y)\|^2 > 
       \sum_{n = 1}^{R^2} \| \phi^{n}_x - \phi^{n}_y \|^2_{p_n} 
  \geq \sum_{n = 1}^{R^2} \| \phi^{n}_x -\phi^{n}_y \|^2_{\infty} = R^2.
\end{equation*}
This concludes the proof.
\end{proof}

\section{Discrete Group Case}
\label{sec:discrete-group-case}

In this section we prove Theorem~\ref{thm:groupcase}.  The proof is
similar to one of the standard proofs that amenable groups are uniformly
embeddable into Hilbert space \cite{bekka-cherix-valette95}.

Let $\Gamma$ be a countable discrete group.  Equip $\Gamma$ with a
proper length function; if $\Gamma$ is finitely generated this could be
a word length function.  The length of an element $s\in\Gamma$ is
denoted $| s |$.  

The left regular representation of $\Gamma$ on $l^p(\Gamma)$ is denoted
$f\longmapsto s\cdot f$; these are isometric actions.  We regard
$l^{\infty}(\Gamma) \subset {\mathcal B}(l^2(\Gamma))$ as multiplication
operators; we regard $\Gamma\subset{\mathcal B}(l^2(\Gamma))$ via the
left regular representation.  With these conventions we have $s f s^* =
s\cdot f$ for all $s \in \Gamma$ and $f \in l^{\infty}(\Gamma)$.

\begin{lemma}
  For every $C>0$ and $\eps > 0$ there exists a non-negative, finitely
  supported function $f$ on $\Gamma$ such that $\|f\|_\infty = f(e) = 1$
  and $\| s\cdot f - f \|_{\infty} < \eps$, for all $s\in\Gamma$
  with $|s|\leq C$.
\end{lemma}

\begin{proof} 
Consider the $C^*$-algebra $B = C^*_r(\Gamma) + 
{\mathcal K}(l^2(\Gamma)) \subset B(l^2(\Gamma))$.  Define 
$h \in c_0(\Gamma)$ by $h(s) = {|s|^{-1}}$ (set $h(e)=1$).  
Observe that $h\in c_0(\Gamma) \subset l^{\infty}(\Gamma) \cap
{\mathcal K}(l^2(\Gamma))$ is a strictly positive element in
${\mathcal K}(l^2(\Gamma))$.  Consequently, $\{\, h^{1/n} \,\}$ is an
approximate unit for ${\mathcal K}(l^2(\Gamma))$ and the convex hull
of $\{\, h^{1/n} \,\}$ contains an approximate unit for 
${\mathcal K}(l^2(\Gamma))$ which is {\em quasicentral} in $B$.  
In particular, given $C>0$ and $\eps>0$ there exists a convex
combination 
%
%
$g$ of the $\{\, h^{1/n} \,\}$ such that 
\begin{equation*}
  \| s\cdot g - g \|_{\infty} = \| s g s^* - g \| 
         = \| s g - g s \| < \frac{\eps}{3},  
\end{equation*}
for all $s \in \Gamma$ with $|s|\leq C$.  Note that $g\in c_0(\Gamma)$
satisfies $\|g\|_\infty = g(e)=1$ and let $f$ be a function of finite
support such that $\| f- g \|_{\infty} < \eps/3$, $\|f\|_\infty= f(e) =
1$.  The proof concludes with a simple estimate:
\begin{equation*}
  \| s\cdot f - f \|_\infty \leq \| s\cdot f - s\cdot g\|_\infty +
                        \| s\cdot g - g \|_\infty + \| g-f \|_\infty
             <\eps,
\end{equation*}
for all $s \in \Gamma$ with $|s|\leq C$.
\end{proof}

\begin{proof}[Proof of Theorem~\ref{thm:groupcase}.]
Apply the previous lemma to a sequence of $C$'s tending to infinity
and a sequence of $\eps$'s tending sufficiently fast to zero to obtain
a sequence of finitely supported functions $\{\, f_n \,\}$ such that
\begin{equation*}
  \sum_{n = 1}^{\infty} \| s\cdot f_n - f_n \|_{\infty}^2 < \infty,
\end{equation*}
for all $s \in \Gamma$. Apply Lemma~\ref{lem:conv} to obtain a
sequence $\{\, p_n \,\}$ such that 
\begin{equation*}
  \sum_{n = 1}^{\infty} \| s\cdot f_n - f_n \|_{p_n}^2 < \infty,
\end{equation*}
for all $s \in \Gamma$.  By virtue of this inequality we may define 
$\Phi : \Gamma \to \oplus l^{p_n} (\Gamma)$ by
$\Phi(s) = \oplus s\cdot f_n - f_n$. 

Denote by $\lambda$ the sum of the left regular representations of
$\Gamma$ on $\oplus l^{p_n}(\Gamma)$; $\lambda$ is an isometric
representation.  A straightforward calculation reveals that $\Phi$
satisfies the cocycle identity with respect to $\lambda$: 
$\Phi(st) = \lambda_s(\Phi(t))+\Phi(s)$
for all $s$, $t\in\Gamma$.  Consequently, $\alpha_s=\lambda_s +\Phi(s)$
defines an affine isometric action of $\Gamma$ on $\oplus
l^{p_n}(\Gamma)$.  Observe that
\begin{equation*}
  \| \Phi(s) -\Phi(t) \| = \| \alpha_s(0) - \alpha_t(0) \| 
               = \| \alpha_{s^{-1}t}(0) \| = \| \Phi(s^{-1}t) \|,
\end{equation*}
for all $s$, $t\in\Gamma$. 

It remains only to check that $\Phi$ is proper, in the sense that
$\|\Phi(s)\|\to\infty$ as $s\to\infty$ in $\Gamma$.
%
This is straightforward.  Indeed,
%
%
if $f_n$ is supported in the ball of radius $R_n$ (in $\Gamma$, w.r.t.\ 
$| \cdot |$) and $s \in \Gamma$ satisfies $|s| > 2R_n$ then the
supports of $s\cdot f_n$ and $f_n$ are disjoint and 
$\| s\cdot f_n - f_n \|_{p_n} \geq 1$ (since $\| f_n \|_{\infty} = 1$). 
Hence, if $|s| > 2max\{R_1,\ldots,R_m\}$ then $\| \Phi(s) \| \geq
m^{1/2}$.
\end{proof}

\end{document}